\documentclass[11pt]{amsart}
\usepackage{amssymb,amsfonts,amsthm,amsmath}
\usepackage{mathrsfs}
\usepackage[english]{babel}

\newtheorem{lemma}{Lemma}[section]

\newtheorem{proposition}[lemma]{Proposition}
\newtheorem{remark}[lemma]{Remark}

\newcommand{\field}[1]{\mathbb{#1}}

\newcommand{\R}{\field{R}}
\newcommand{\Sr}{\field{S}}

\newcommand{\Rn}{\R^n}

\newcommand{\nf}{\nabla f}

\newcommand{\barvp}{\overline{\varphi}}
\newcommand{\bB}{{\bf B}}
\newcommand{\bg}{{\bf g}}
\newcommand{\bh}{{\bf h}}

\newcommand{\bu}{{\bf u}}
\newcommand{\bx}{{\bf x}}
\newcommand{\clos}{{\rm{\bf clos}}}

\newcommand{\oo}{{\bf 0}}
\newcommand{\reg}{{\rm reg}}
\newcommand{\tildbeta}{{\widetilde{\beta}}}
\newcommand{\ud}{{\rm d}}
\newcommand{\vp}{{\varphi}}
\numberwithin{equation}{section}

\begin{document}
\title[singular gradient trajectories and oscillations]
{gradient trajectories for plane singular metrics I: oscillating trajectories}

%    Information for author
\author{Vincent Grandjean}
\address{{\it Temporary Address:} Departamento de Matem\'atica, UFC,
Av. Humberto Monte s/n, Campus do Pici Bloco 914,
CEP 60.455-760, Fortaleza-CE, Brasil}
\email{vgrandje@fields.utoronto.ca}
%\author{Fernando Sanz}
%\address{}
%\email{}

\thanks{The author is very grateful to K. Bekka and L. Paunescu for having
independently brought this question to his attention, and to A. Fernandes for 
remarks and comments. The author is pleased to thank IRMAR of Universit\'e de Rennes I, 
for support while working on this note.}
% contract number HPRN-CT-2001-00271}
% \thanks{Partially Supported by Deutsche Forschungs-Gemeinschaft in the Priority Program
% \it Global Differential Geometry}

%    General info

% \date{\today}

\begin{abstract}
In this short note, we construct an example of a real plane analytic singular metric, 
degenerating only at the origin, such that any gradient trajectory (respectively to this singular metric) of 
some well chosen function spirals around the origin. The inversion mapping carries this example into an 
example of a gradient spiraling dynamics at infinity.
\end{abstract}
\maketitle
%
%
%
%
%
%
%     *****************************************************************
%
%
%
%
%
%
%.
%
%
\section{Introduction}
In the early 60s, Thom asked about the behaviour of the (Euclidean) gradient flow of 
a given real analytic function nearby the critical locus of the function. He conjectured that 
any gradient trajectory with limit point a critical point $\oo$ (the origin) 
should have a limit of secants at the origin. It took around thirty years to 
eventually prove that \em Thom Gradient Conjecture \em was true. This was achieved by Kurdyka, Mostowski 
and Parusi\'nski \cite{KMP}, using intensively \L ojasiewicz's result on the finiteness of the 
length of gradient trajectories in a neighbourhood of a limit point \cite{Loj1}. 
Nowadays questions around the dynamics of a gradient trajectory or of a pencil of trajectories nearby 
a limit point have switched to asking whether they are analytically oscillating or not.
A gradient trajectory is analytically \em non-oscillating \em if for any semi-analytic subset,
the intersection with the given gradient trajectory has finitely many connected components. 
Given any real analytic gradient differential equation, Moussu Theorem \cite{Mou} ensures that 
there are always regular analytic curves through a singular point such that each open half-branch
of the curve is a gradient trajectory. Such trajectories are obviously non-oscillating. 
But a few other particular cases \cite{Sa,FS,Go,GS} the presence of other non-oscillating trajectories 
is not known.   

\smallskip
Assume that $(\Rn,\oo)$ is equipped with any given real analytic Riemannian metric $\bg$.
Given a real analytic isolated surface singularity germ $(S,\oo)$ of $(\Rn,\oo)$, the 
regular part $S_\reg$ of $S$ is equipped with $\bg_S$ the restriction of the ambient Riemannian metric
to the surface. For any real analytic function $f:(\Rn,\oo) \to (\R,0)$, possibly singular at $\oo$, 
we can consider its restriction $f_S$ to the surface $S$, and consider the gradient vector field $\nf_S$ of the
function $f_S$ (relative to the metric $\bg_S$). It is a real analytic vector field, defined only on $S_\reg$.
A trajectory of the restricted gradient $\nf_S$ will be called a \em restricted gradient 
trajectory. \em The presence of the singular point $\oo$ at which the metric $\bg_S$ cannot be extended 
may \`a-priori considerably influence the dynamics of the restricted trajectories nearby $\oo$.
Nevertheless, this problem is completely understood in the joint work of the author and F.Sanz \cite{GS}. 
Our result states that restricted gradient trajectories do not oscillate at their limit point $\oo$.
\\
Since \L ojasiewicz \cite{Loj2} we know that the germ $(S,\oo)$ is topologically a finite 
union of the closure of positive cones over a circle with vertex $\oo$ and so the germ $((S\setminus \oo),\oo)$ 
has finitely many connected components, each of which is a positive cone. 
In particular, close enough to the singular point $\oo$ any restricted gradient trajectory of $\nf_S$ 
must stay in a single positive cone. Such a positive cone is just a punctured plane. 
From a purely topological dynamical point of view, the problem of the local behaviour of a 
restricted gradient trajectory nearby the singular point looks like a particular case of the local 
behaviour of the trajectories of a vector field which outside the point $\oo$ is 
the gradient vector field of a real analytic Riemannian metric $\bh$, 
outside the origin $\oo$. If the metric $\bh$ cannot be extended through $\oo$
into a Riemannian metric, then the metric $\bh$ is not positive definite at $\oo$ and
we will call $\oo$ the \em degeneracy locus \em of $\bh$.   
A trajectory of such a vector field is a \em singular gradient trajectory. \em
Although there are similarities with the result of the author joint work \cite{GS}, 
the local behaviour of trajectories of singular gradient differential equations 
at a point of the degeneracy locus of the singular metric, is far wilder 
than that of restricted gradient trajectories as we will see.

The question of the oscillation of singular gradient trajectories 
was asked to us simultaneously by Prof. K. Bekka and Prof. L. Paunescu.

\smallskip
The paper is organized as follows: 

Section 2 introduces in a wider context, the problem of studying the (singular) gradient trajectories 
of a function relative to a singular metric nearby the degeneracy locus of the metric, 
that is the locus of points where the corresponding $2$-symmetric tensor is not positive definite 
is not empty. 

Sections 3 and 4 are devoted to build an example of singular gradient of a function relative 
to a singular metric degenerating only at the origin, such that the corresponding (singular) 
gradient trajectories spiral, thus oscillate, in a neighbourhood of a point of the degeneracy locus. 
We proceed along the following lines:\\
We will build a real analytic Riemannian metric $\bh$ onto the punctured unit ball $\bB_1^*:= \bB_1  
\setminus \oo$ which extends into a real analytic $2$-symmetric tensor through the origin.
To achieve that, we first build a real analytic $2$-symmetric tensor $\bg$
on the spherical blowing-up $[\bB_1,\oo]$ of the disk $\bB_1$ such that
it is a Riemannian metric onto the pull-back of the punctured disk and 
is only positive semi-definite along the boundary circle, exceptional locus of the blowing-up.
Thus we find a real analytic function on $[\bB_1,\oo]$ whose gradient
trajectory accumulates along the whole boundary circle.
Up to a rescaling of $\bg$ by a non-negative function vanishing only on
the boundary circle, we blow-down this singular metric to find the singular metric 
$\bh$, and then we blow-down the function. Both are real analytic. Since (singular)
gradient differential equation are only sensitive to the conformal structure 
of the (singular) metric, we are guaranteed that the singular gradient trajectory of the 
blown-down function spiral around the origin (Proposition \ref{prop:spiraling}).

Section 5, although short, exploits the previous counter-example using the inversion mapping, 
to exhibit a plane smooth semi-algebraic metric in a neighborhood of infinity for which 
there exists a smooth semi-algebraic function with a spiraling gradient dynamics at infinity 
(Proposition \ref{prop:spiral-infinity}).

In the last section we speculate about which properties of the metric at the singular 
point could cause the oscillating phenomenon, when the geometry of the function 
is too special in regards of that of the singular metric. 
%
%
%
%
%
%
%
%
%
%
%
%
%     *****************************************************************
%
%
%
%
%
%
%
%
%
%
%
%
\section{On singular gradient differential equations}
Let $M$ be a real analytic connected manifold. 
A real analytic $2$-symmetric tensor $\bh$ defined on $M$ is called a \em singular metric \em
on $M$, if there exists a real analytic subset $Y$ of $M$ of codimension larger or equal to $1$, 
such that $\bh$ is positive definite on $M\setminus Y$, and degenerates along $Y$, that is at each point $y\in Y$, 
the quadratic form $\bh(y)$ is only semi-positive definite. The subset $Y$ is called
\em the degeneracy locus of the singular metric \em $\bh$.

\medskip
Given a real analytic function $f:M\to \R$, we consider the vector field $\nabla_\bh f$ defined on $M\setminus Y$ as 
dual of the differential $\ud f$ for a given singular metric $\bh|_{M\setminus Y}$. 
By definition, we obtain
\vspace{4pt} 
\begin{equation}
\vspace{4pt}
\ud_x f \cdot u = \langle \nabla_\bh f(x), u\rangle_\bh, \; \forall x \in M\setminus Y, \; \forall u \in T_x M,
\end{equation}
where $ \langle, \rangle_\bh$ denotes the scalar product coming from $\bh$.
Once are given some coordinates $x$ nearby a point $\bx_0 \in M$, the quadratic form $\bh(x)$ is given 
by a matrix $H(x)$. The vector field $\nabla_\bh f$ is given in the local coordinates by $H^{-1}(x) \partial f (x)$, 
where $\partial f$ is the vector fields of the partial derivative of $f$ in the local coordinates.
Let $H^*$ be the adjoint matrix of the matrix $H$. We recall that $H(x) \cdot H(x)^* = H(x)^* \cdot H(x)  =
\det H(x) Id = \det \bh (x) Id$. In local coordinates we can define the vector field 
\vspace{4pt} 
\begin{equation}
\vspace{4pt}
\xi_\bh f (x):= H^*(x) \cdot \partial  f (x) = \det \bh (x) \nabla_\bh f (x)
\end{equation} 
The vector field $\xi_\bh f$ is independent of the coordinates chosen.
It is a real analytic vector field on the whole of $M$, 
co-linear to $\nabla_\bh f$ and vanishing on the subset $Y$. 

With an obvious abuse of language we call the next 
differential equation the \em singular gradient differential equation of the function $f$ relative to
the singular metric $\bh$: \em
\vspace{4pt} 
\begin{equation}\label{eq:degen-grad}
\vspace{4pt}
\dot{x} (t) = \xi_\bh f (x(t)), \: x(0) = x_0 \notin Y.
\end{equation}

We would like to inquire about the behaviour of the singular gradient trajectories
in a neighbourhood of a point $y$ of $Y$. 

In the present paper we are going to provide a simple example of such a situation where 
the degeneracy locus $Y$ is the origin of $M$ the real plane, and such that all the trajectories accumulate
at this point in oscillating.
%
%
%
%
%
%
%
%
%
%     *****************************************************************
%
%
%
%
%
%
%
%

\section{Plane counter-example:how to make it}\label{Sec:HTMI}
We will provide an example of a singular metric $\bh$ on the real plane whose degeneracy locus 
is just the origin, and find a function for which all singular gradient trajectories 
spiral around the origin.
\\
In order to do so, we will work on a spherical blowing-up of the plane. We will produce there
a singular metric $\bg$ degenerating only on the boundary circle, that is the pre-image 
of the degeneracy locus of $\bh$ under the spherical blowing-up.
The singular metric $\bg$ we will choose will be, up to the multiplication by a function vanishing 
only along the boundary circle, the pull-back under the spherical blowing-up
of the singular metric $\bh$ we want to produce. 
The rescaling factor of the lifting of our singular metric $\bh$ is of no importance since 
what the foliation induced by a (singular) gradient differential equation uses from the 
(singular) metric is its conformal structure. Such a gradient foliation is insensitive to the 
sole change of the measure of the length.

\bigskip
Let us consider the spherical blowing up of the real plane: 
\begin{center} 
\vspace{4pt}
$\beta : \Sr^1\times \R_{\geq 0}\to \R^2$, defined as $(\bu,r) \to r\bu$.
\vspace{4pt}
\end{center}
The pre-image $\beta^{-1}(\oo)$ of the origin $\oo$ is thus the boundary circle $\Sr^1 \times 0$.
Instead of working exactly on $\Sr^1\times\R_{\geq 0}$, we are going to work on its 
universal covering $\R\times\R_{\geq 0}$ to exhibit the metric and see very well the 
spiraling behaviour around the boundary circle of the singular gradient trajectories 
relatively to the singular metric we will consider. Thus the boundary circle $\Sr^1\times 0$ 
is replaced by the boundary line $\R\times 0$.

We consider the following $2$-symmetric tensor on $\R\times [0,1[$:
\vspace{4pt}
\begin{equation}\label{eq:degen-metric}
\vspace{4pt}
\bg = \ud r^2 + 2r^3\ud r\ud \vp + r^4\ud \vp^2,
\end{equation}
in the coordinates $(\vp,r)$ in $\R\times [0,1[$.
\\
Given any $(u,v) \in \R^2 \setminus \oo$, we check easily that for any $r\in ]0,1[$
the real number $u^2 +  2r^3 u v + r^4 v^2$ is positive.
\\
The determinant of this metric is $r^4(1-r^2)$ and thus vanishes on $r=0$. Thus $\bg$ is a 
Riemannian metric on $\{0<r<1\}$ and degenerate along the boundary line.

\medskip
Given any smooth function $(\vp,r) \to f(\vp,r)$ defined over $\R\times \R_{\geq 0}$ the gradient vector 
field $\nabla_\bg f$ of $f$ for the degenerate metric $\bg$ is 
\begin{center}
\vspace{4pt}
$r^4(1-r^2) \nabla_\bg f = [r^4\partial_r f - r^3 \partial_\vp f]\partial_r + [-r^3\partial_r f + 
\partial_\vp f]\partial_\vp$.
\vspace{4pt}
\end{center}
Thus the gradient differential equation associated with $f$, up to multiplication by $r^4$ reads 
\vspace{4pt}
\begin{equation}
\vspace{4pt}
\left\{
\begin{array}{rcl}
\vspace{6pt}
\dot{r} & = & r^4\partial_r f - r^3 \partial_\vp f \\
\dot{\vp} & = & -r^3\partial_r f + \partial_\vp f
\end{array}
\right.
\end{equation}  

Let us see how the solution of this differential equation does behave nearby $\{r=0\}$ in the very simple case 
of $f (\vp,r)= -r$. \\
Namely it reduces to  
\vspace{4pt}
\begin{equation}
\vspace{4pt}
\left\{
\begin{array}{rcl}
\vspace{6pt}
\dot{r} & = & -r \\
\dot{\vp} & = & 1
\end{array}
\right.
\end{equation}    
We deduce that any trajectory from a point $(\vp_0,r_0)$ with $r_0>0$ never ends-up on a point $(\vp_1,0)$
since $r=0$ is a trajectory of the above differential equation.\\
We can check  that for $r>0$ this differential equation reads
\vspace{4pt}
\begin{equation}
\vspace{4pt}
\displaystyle{\frac{\ud \vp}{\ud r} = -\frac{1}{r}},
\vspace{4pt}
\end{equation}
Thus any trajectory from a point $(\vp_0,r_0)$ is a graph in $r$ of a function $\vp(r) = C_0 -\ln(r)$ and thus, 
$\vp (r)$ tends to $+\infty$ as $r$ tends to $0$.
\begin{remark}\label{rk:further-examples}
In the present case we observe that any function $g$ of the form $g:=f + r^4 h$, for any real analytic 
function $h$ defined in a neighbourhood of the boundary circle, will provide singular gradient
trajectories which accumulate at a point of $\R \times 0$ only at infinity.
\end{remark}
%
%
%
%
%
%
%
%
%
%
%
%
%      **********************************************************
%
%
%
%
%
%
%
%
%
%
%
%
\section{Plane singular gradient trajectories spiraling around the origin}
As described at the beginning of Section \ref{Sec:HTMI} we have explained how to provide 
the singular metric on the plane. We will use the blowing-down mapping $\beta$, in polar coordinates, 
in order to find the singular metric that will give the spiraling-around-the-origin behaviour of the 
whole phase portrait of the singular gradient trajectories of the function Euclidean distance to the origin 
for the singular metric $\bh$.

\medskip
Let $(x,y)$ be coordinates in $\R^2$, and let us write $x = r\cos \barvp$ and $y = r \sin \barvp$,
For $(r,\barvp) \in \R_{\geq 0} \times [0,2\pi]$.

We thus find 
\begin{center}
$\begin{array}{rcl}
\vspace{6pt}
r\ud r & = & x \ud x + y\ud y \\
r^2\ud \barvp & = & x\ud y - y\ud x
\end{array}
$
\vspace{4pt}
\end{center}  
Defining $\bh (x,y)$ as $r^2 \bg (r,\barvp)$ we find 
\begin{center}
$\begin{array}{rcl}
\vspace{6pt}
\bh & = & (x\ud x + y\ud y)^2 + 2 r^2(x\ud x + y\ud y)(x\ud y - y \ud x) + r^2(x\ud y - y \ud x)^2 \\ 
\vspace{6pt}
& = & [x^2 + r^2(-2xy + y^2)]\ud x^2 + 
2[xy+r^2(x^2-y^2-2xy)]\ud x \ud y + \\
& & 
[y^2 + r^2(2xy + x^2)]\ud y^2
\end{array}
$
\vspace{4pt}
\end{center} 
Thus $\bh$ defines a real analytic $2$-symmetric tensor on $\R^2$ whose degeneracy locus is the origin at 
which it is the null quadratic form. Note that $\bh$ is positive definite on $\R^2 \cap \{0<r<1\}$.
Consequently we restrict our attention to the open unit ball $\bB_1$ of $\R^2$,
which we equip with the singular Riemannian metric $\bh$ just defined above.

\medskip
Let us consider the universal covering of the spherical blowing-up of $\bB_1$, namely, 
$\tildbeta :[0,1[\times \R \to \bB_1$, defined as $(r,\vp) \to (r\cos \vp,r \sin \vp)$. \\
Thus for any interval $I = [a,a+2\pi] \subset \R$, for any real number $a$, the restriction of $\tildbeta$ to 
$[0,1[ \times I$ induces a diffeomorphism $]0,1[\times I$ onto the punctured ball $\bB_1^*$.\\
We obviously check that $\tildbeta^* (\bh)  = r^2 \bg$.
We want to understand the asymptotic behaviour of the gradient differential 
equation $\dot{p} = \nabla_\bh \delta (p)$, defined on $\bB_1^*$, nearby the boundary of this 
domain, namely the origin $\oo$. In the coordinates $(x,y)$ this differential equation reads as 
\begin{center}
\vspace{4pt}
$\left\{
\begin{array}{rcl}
\vspace{6pt}
\dot{x} & = & -2x[y^2 + r^2(2xy + x^2)]+2y[xy+r^2(x^2-y^2-2xy)] \\
\dot{y} & = & 2x[xy+r^2(x^2-y^2-2xy)] - 2y[x^2 + r^2(-2xy + y^2)]
\end{array}
\right.
$
\vspace{4pt}
\end{center}

When pulled back by $\tildbeta$ this differential equation transforms into the differential equation of the gradient of 
$(r,\vp) \to r^2$, that is the differential equation given by the vector field $2r\nabla_\bg r$. Thus its trajectories are 
the same as that of  $\nabla_\bg r$ in $\{r>0\}$.
\\
Moreover, any non stationary trajectory of $\nabla_\bh \delta$ is lifted by $\tildbeta$ in a unique trajectory of 
$\nabla_\bg r$ lying in $\{r>0\}$, the converse is also true. 
\\
But as we have already checked in the third Section, any trajectory of $\nabla_\bg r$ with initial data
lying in $r>0$ is a curve of the form $r \to (C_0-\ln r, r)$. Thus the image of such gradient trajectory 
will be mapped by $\tildbeta$ on a gradient trajectory of $\nabla_\bh \delta$ lying in 
$\bB_1^*$ and parameterized as $r \to (r \cos (C_0 - \ln r), r \sin (C_0 - \ln r))$ which spirals 
around the origin as $r$ goes to $0$. Thus we have proved the following 
\begin{proposition}\label{prop:spiraling}
Any singular gradient trajectory (respectively to the singular Riemannian metric $\bh$) of the function 
\begin{center}
\vspace{4pt}
$\delta:\bB_1 \to \R$ defined as $(x,y) \to \delta (x,y) := - (x^2+y^2) = - r^2$, 
\vspace{4pt}
\end{center}
with initial data $(x_0,y_0) \neq \oo$, spirals around the 
origin as the "time" goes to $+\infty$.
\end{proposition}
\begin{remark}\label{rk:further-spiraling}
To echo Remark \ref{rk:further-examples}, any function $g$, analytic or not, 
of the form $-r^2 + r^5h$, for $h$ a $C^1$ function in a neighbourhood of $\oo$, will 
provide singular gradient trajectories, relative to the singular metric $\bg$
above, which spiral around the origin $\oo$.
\end{remark}
%
%
%
%
%
%
%
%
%           ********************************************************
%
%
%
%
%
%
%
%
%
%
\section{Example of a spiraling gradient dynamics at infinity}
Let $f:\Rn\to \R$ be a $C^2$ semi-algebraic function. The main result of \cite{Gr2} is that any Euclidean gradient 
trajectory of the function $f$ leaving any compact subset of $\Rn$ has a limit of secants at infinity.

\medskip
If the space $\Rn$ is equipped with a Riemannian analytic metric 
the behaviour of a half-gradient curve is either to accumulate on a point in $\Rn$ or
to leave any compact subset of $\Rn$. 
The image of $\Rn$ under the smooth semi-algebraic diffeomorphism $p\to\frac{p}{\sqrt{1+|p|^2}}$ is the 
open Euclidean unit ball $\bB_1$. The sphere bounding this unit ball will be referred to as \em
the sphere at infinity. \em
The behaviour of the given Riemannian metric nearby the sphere at infinity is of great importance 
for the respective gradient curve leaving any compact subset. 
\\
If we restrict our attention to the plane, intuition and reason command that there should exist
many metrics in the neighbourhood of the circle at infinity for which we should find 
globally subanalytic analytic functions whose gradient trajectories leave any compact and spiral, 
in other words accumulates at each point of the boundary circle at infinity. We are going to give 
such an example below.

\medskip
Let $\bh$ be the singular metric  of Section 4 defined on $\bB_1$. 
The smooth semi-algebraic function $-f:\R^2\to \R$, where $f(p) := \frac{|p|^2}{1+|p|^2}$, has  
all its singular gradient trajectories respectively to $\bh$ which accumulate onto 
the origin in spiraling. In other words given any singular gradient trajectory $\gamma$, 
for any unit vector $\bu\in \Sr^1$, there exists a sequence $(p_k)_k$ of points of $\gamma$ 
such that $p_k\to \oo$ and $\frac{p_k}{|p_k|} \to \bu$ as $k\to\infty$.   

\smallskip
Let us consider the plane inversion mapping:
\begin{center} 
\vspace{4pt}
$I :\R^2\setminus \oo \to \R^2\setminus\oo$ defined as $p\to \frac{p}{|p|^2}$.
\vspace{4pt}
\end{center}

We find that $I^*\bh$ is a smooth semi-algebraic Riemannian metric on
the pre-image $I^{-1}(\bB_1) = \R^2\setminus\clos(\bB_1)$. We find that $I^* f =f$ and also observe 
that the origin, seen as the boundary circle of the punctured unit disk, is "mapped" onto the boundary 
circle at infinity. Thus we deduce the following
\begin{proposition}\label{prop:spiral-infinity}
Any gradient curve $\gamma$ of the function $I^*f$ for the metric $I^*\bh$ leaves any compact subset of $\R^2$ and 
acccumulates on the whole boundary circle at infinity, in other words for any unit vector $\bu\in \Sr^1$, there exists 
a sequence $(p_k)_k$ of points of $\gamma$ such that $|p_k|\to \infty$ and $\frac{p_k}{|p_k|} \to \bu$ as $k\to\infty$.    
\end{proposition} 
%
%
%
%
%
%
%
%
%           ********************************************************
%
%
%
%
%
%
%
%
%
%
\section{Remarks and Speculations}
1) Given a plane real analytic singular metric, there will be nevertheless uncountably many 
functions whose singular gradient trajectories will not spiral, for instance
those taking positive and negative values close to the origin.
Now given a real analytic function vanishing only at the origin, finding necessary and sufficient conditions 
so that the pair singular metric and function does not rise a singular gradient differential 
equation with a spiraling dynamics around $\oo$ seems for the moment complicated.
If we are concerned only by the properties of the metric this problem is partially solved in
\cite{Gr2}.

\bigskip
2) The topological description of the author joint work \cite{GS} of restricted gradient on isolated
surface singularity suggested that a singular metric degenerating only at a single point 
might have produced non-oscillating singular gradient trajectories. As we show here it is generally not true. 
In particular this naive point of view is forgetting that 
in this description the degeneracy of the metric is forced by the space, since the restricted metric can 
extend to the whole ambient space as a standard Riemannian metric. Consequently
The singular metric of \cite{GS} comes from the restriction of a Riemannian metric to
a singular "cone". Note that the asymptotic behaviour at the singular point of this restricted metric
is just the restriction of the ambient metric to the asymptotic behaviour of the singular "cone" at 
its tip, namely the limits at the tip of the tangent spaces to the surface.

\bigskip
3) The singular metric we have exhibited here presents two asymptotic behaviours at $\oo$ which may play 
some role in the spiraling example presented. First any limit at $\oo$ of the normalized quadratic forms 
associated to the metric are of rank $1$. Consequently to the remark of point 2), the example presented here 
clearly forbid the embedding of $B_1$ in $\R^n$ so that the metric $\bh$ is the restriction of an ambient metric.
Moreover in the projective space of quadratic forms the sets of these limits is an embedded projective line.
In other words, the "conformal" structure carried by the singular metric $\bh$, which is the 
quality that matters for (singular) gradient trajectories, cannot be defined uniquely at the origin 
since there is an embedded projective line of such possible limits.
Then some geometric properties of the function are such that its singular gradient trajectories
spiral around the origin. Which are they relatively to the singular metric, we just do not know yet.
%
%
%
%
%
%
%
%
%
%
%
%
%
%      **********************************************************
%
%
%
%
%
%
%
%
%
%
%
%

\end{document}